\documentclass{amsart}

\def\bd{{\bf d}}

\usepackage{graphicx}
\usepackage{color}
\newtheorem{theorem}{Theorem}[section]
 \newtheorem{lemma}{Lemma}[section]
 \newtheorem{prop}{Proposition}[section]
 \newtheorem{defn}{Definition}
  \newtheorem{remark}{Remark}[section]

\UseRawInputEncoding
\usepackage[russian, english]{babel}

\def\br{{\bf r}}
\def\dfrac#1#2{\displaystyle{#1\over #2}}

\def\bv{{\bf V}}
\def\bV{{\bf V}}

\def\Div{\mbox{div}\,}
\def\Rot{\mbox{curl}\,}

\def\bB{{\bf B}}

\def\bE{{\bf E}}

\begin{document}


%
%

\title[A 
Electrostatic oscillations of an electron plasma]{
On the properties of multidimensional electrostatic oscillations of
an electron plasma
}

\author{Olga S. Rozanova}

\address{ Mathematics and Mechanics Department, Lomonosov Moscow State University, Leninskie Gory,
Moscow, 119991,
Russian Federation,
rozanova@mech.math.msu.su}



\maketitle


\begin{abstract}
We consider the classical Cauchy problem for a system of equations
describing 3D arbitrary electrostatic oscillations of the cold
plasma and introduce an iteration procedure that allows estimating the blow-up time from below.
This procedure is constructive provided one succeed to obtain a two-sided estimate of an additional quantity depending on the solution.
We show that this is possible in the case of one and two dimensions, as well as for solutions with zero vorticity.
 For the particular case of two-dimensional initial data with the radial symmetry, refined sufficient
 conditions for destruction and preservation of smoothness in the first  period of oscillations are obtained. Moreover, we give the example of estimating the blow-up time for the data for which results of numerics exist and discuss a roughness of our estimate.

 \end{abstract}

\section{Introduction}


The equations of hydrodynamics of "cold" or electron plasma in the
non-relativistic approximation in dimensionless quantities take the
form(see, e.g.,~\cite{ABR78},  \cite {GR75})
\begin{eqnarray}\label{base1.1}
\dfrac{\partial n }{\partial t} + \Div(n \bv)=0,\quad
\dfrac{\partial \bv }{\partial t} + \left( \bv \cdot \nabla \right) \bv
= \, - \bE -  \left[\bv \times  \bB\right],\\
\frac{\partial \bE }{\partial t} =   n \bv
 + {\rm rot}\, \bB,\qquad 
\frac{\partial \bB }{\partial t}  =
 - {\rm rot}\, \bE,\qquad \Div \bB=0,\label{base1.5}
\end{eqnarray}
$ n$ and $ \bv=(V_1, V_2, V_3)$  are the density and velocity of electrons,
$ \bE=(E_1, E_2, E_3)$ and  $  \bB=(B_1, B_2, B_3) $  are vectors of electric and magnetic fields.
All components of solution depends on $t\in {\mathbb R}_+$ and $x\in{\mathbb R}^3$.

At present, much attention is paid to the study of cold plasma in
connection with the possibility of accelerating electrons in the
wake wave of a powerful laser pulse \cite{esarey09}; nevertheless,
there are very few theoretical results in this area.

It is commonly known that the plasma oscillations described by \eqref{base1.1}, \eqref{base1.5}, tend to break. Mathematically,
 the breaking process means a blow-up  of the solution, and the appearance of a delta-shape singularity of the electron density  \cite{Dav72}. Among main interests  is a study of possibility of existence of a smooth solution as long as possible.

Let us assume that the oscillations  are electrostatic, i.e. ${\rm
rot}\, \bE$ and the magnetic field $\bB$ does not change with time.
For simplicity, we set  $\bB=0$.

 From the first equations of \eqref{base1.1} and  \eqref{base1.5} under the assumption
 that the solution is sufficiently smooth and that the steady-state density is equal to 1, it follows
\begin{eqnarray}n=1- \Div \bE,\label{n}\end{eqnarray}
thus, $ n $ can be removed from the system. Thus, we get
\begin{eqnarray}\label{4}
\dfrac{\partial \bv }{\partial t} + \left( \bv \cdot \nabla \right)
\bv = \, - \bE,\quad \frac{\partial \bE }{\partial t} + \bv \Div \bE
 = \bV,
\end{eqnarray}
together with
\begin{equation}\label{cond}
\Rot \bE=0,\quad  \Rot ((1-\Div \bE)\bv)=0.
\end{equation}

Consider the initial data
\begin{equation}\label{CD1}
(\bv, \bE) |_{t=0}=
(\bv_0, \bE_0)(x)  \in C^2({\mathbb R}^3),\quad x\in {\mathbb R}^3. 
\end{equation}
with properties \eqref{cond}.

\begin{defn}
We will say that {\it the solution to the problem \eqref{4}, \eqref{CD1} does not blow up}
(or the oscillations do not break) on $[0,T),\, T>0,$ if the density $n$, found as \eqref{n} remains bounded for all $t\in T.$
\end{defn}

The main problem is that the hyperbolic system \eqref{4},
\eqref{CD1} has locally in time a smooth solution \cite{Daf16},
however, it not necessary satisfies \eqref{cond}. In the general
case system \eqref{4}, \eqref{cond} is overdetermined.

Nevertheless, as it will be shown below, there exist important
classes of solutions such that \eqref{cond} automatically hold.
It is affine solutions and radially symmetric solutions (see Sec.\ref{Sec_aff_ax}). 
Therefore below we {\it assume } that we deal only with solutions
with property  \eqref{cond}.

We denote ${\mathcal D}=\Div \bv$, ${\bf \Xi}=(\xi_1, \xi_2,
\xi_3)=\Rot \bv$, $\lambda=\Div \bE$,
$J_{1}=\det(\|\partial_{x_i}{V_j}\|),$ $i,j=1,2$,
$J_{2}=\det(\|\partial_{x_i}{V_j}\|),$ $i,j=2,3$,
$J_{3}=\det(\|\partial_{x_i}{V_j}\|),$ $i,j=1,3$, $J=J_1+J_2+J_3.$

System \eqref{4},  \eqref{cond} implies
\begin{eqnarray*}
\dfrac{\partial {\mathcal D}}{\partial t} &+& ( \bv \cdot \nabla
{\mathcal D} ) =   - {\mathcal D}^2 + 2J -\lambda,\qquad
\dfrac{\partial \lambda }{\partial t} + ( \bv \cdot \nabla  \lambda ) ={\mathcal D}(1-\lambda),\\
\end{eqnarray*}
 Thus, along the characteristics $\dfrac{dx_i}{dt}=V_i,\, i=1,2,$ starting
 from the point $x_0$ we obtain the Cauchy problem for the nonlinear system of three ODEs  (non-closed due to $J$):
\begin{eqnarray}\label{base4.3}
 \dot{\mathcal D } =   - {\mathcal D}^2 + 2J-\lambda,\quad
\dot \lambda ={\mathcal D}(1-\lambda),\\
({\mathcal D }, \lambda(t, x_0)|_{t=0}=({\mathcal D}_0, \lambda_0(x_0).\nonumber
\end{eqnarray}

\begin{defn}
We will say that {\it the breaking of oscillations does not occur in
the first period} for the initial data \eqref{CD1}, if the
projection of each phase trajectory starting at
 $(\lambda_0, {\mathcal
D}_0)= $ $ ({\rm div} \bE _0(x_0), {\rm div} \bv_0(x_0))$, $x_0\in
{\mathbb R}^3$ on the plane $(\lambda, \mathcal D)$, intersects the
semi-axis ${\mathcal D}=0$, $\lambda<0$ at least at one point (the
starting point of the trajectory does not count).
\end{defn}
In other words, all trajectories make at least one revolution around
the coordinate origin on the plane $(\lambda, \mathcal D)$.

In Sec.\ref{Sec2} we show that if a two-sided estimate
\begin{equation*}
J_-( \lambda, \mathcal D^2) \le  J\le J_+( \lambda, \mathcal D^2)
\end{equation*}
is known, it is possible to construct an iterative procedure that allows one to find the guaranteed number of oscillations before the blow-up (Sec.\ref{iter}) and to estimate from below the lifetime of a solution with bounded density (Sec.\ref{time}).
In
Sec.\ref{Sec3} we consider particular cases of the problem \eqref{4},
\eqref{CD1} and show that in some cases the solution remains smooth
for all $t>0$.
In
Sec.\ref{Sec4}  is dedicated to the cases, where
the estimate of $J$ can be obtained. For  the case of plain oscillations (Sec.\ref{Sec4.1}) and the case of irrotational oscillations  in any dimensions
(Sec.\ref{Sec4.2}), we prove a sufficient condition for the boundedness of the density in the first period of oscillations.
 In Sec.\ref{Sec5} we construct two-sided estimates of the projection of the integral trajectory to the plane $(\lambda, {\mathcal D})$ for the case of radially symmetric oscillations in any dimensions.
In Sec.\ref{Sec6} we consider the initial data in the
 form of a standard laser pulse corresponding to the 2D (plain) axisymmetric solution and prove refined estimates allowing us to obtain, in particular,  sufficient conditions for the preservation of smoothness by the solution to the Cauchy problem on several periods of oscillations and compare it with the results of direct numerics.
In Sec.\ref {SecD}, we summarize the results obtained and discuss
further hypotheses and open problems.

\section{Method for obtaining a sufficient conditions for boundedness of density}\label{Sec2}

 We denote for convenience
 $ s = \lambda-1 \le 0 $ (the inequality follows from \eqref {n}, since $ n \ge 0 $). Thus, we get
\begin{eqnarray*}
\frac{d\mathcal D}{d s}=\frac{{\mathcal D}^2 +2s  +2J+ 1}{{\mathcal D}s}.
\end{eqnarray*}
Further, replacing  $Z={\mathcal D}^2\ge 0$, we obtain
\begin{eqnarray}\label{ZJ}
\frac{d Z}{d s}=\frac{2 (Z +s + 1)}{s} +\frac{4J}{s}\equiv Q(s,Z,J).
\end{eqnarray}
Let us {\it assume} that the following inequality is known:
\begin{eqnarray}\label{J_est}
 2 J_-(Z,s) \le  2J \le 2 J_+(Z,s).
\end{eqnarray}
Then, taking into account the sign of $s$, we obtain the estimate:
\begin{eqnarray}\label{eqtQ}
 Q_1(s,Z) \le Q(s,Z,J)\le  Q_2(s,Z),
\end{eqnarray}
 where
\begin{eqnarray}\label{Qm}
 Q_1(s,Z) = \frac{2 (Z +s + 2 J_+(s,Z) +1)}{s},
\end{eqnarray}
\begin{eqnarray}\label{Qp}
 Q_2(s,Z) = \frac{2 (Z +s + 2 J_-(s,Z) +1)}{s}.
\end{eqnarray}

Now we can apply Chaplygin's theorem on differential inequalities,
according to which the solution $ Z (s) $ of the Cauchy problem for
\eqref {ZJ} with initial conditions $ Z (s_0) = Z_0 $ for $ s >s_0 $
satisfies the inequality
\begin{eqnarray}\label{Z_ineq1}
  Z_1(s) \le Z(s, J)\le  Z_2(s),
\end{eqnarray}
and for $s<s_0$ the inverse inequality
\begin{eqnarray}\label{Z_ineq2}
  Z_2(s) \le Z(s, J)\le  Z_1(s),
\end{eqnarray}
where $Z_k (s)$ are the solutions to problems $ \frac{d Z}{d s}=Q_k (s,Z)$, $Z(s_0)=Z_0 $, $k=1,2.$

 Since $Z\ge 0,$ then the phase trajectory outgoing from the point  $(s_0, Z_0) $, is bounded,
 if the solution $Z_1(s)\to -\infty$  as  $s\to -\infty$.
If $Z_2(s)\to +\infty$  as  $s\to -\infty$, without becoming zero,
then the phase trajectory outgoing from this point is unbounded.

Thus, to obtain sufficient conditions for the global in time boundedness
of the density and sufficient blow-up conditions, it is enough for us to solve
the equations $ \frac{d Z}{d s}=Q_1 (s,Z)$ and $ \frac{d Z}{d s}=Q_2 (s,Z)$.

First of all, note that on the plane $ (s, {\mathcal D}) $ the
motion is clockwise (for positive $ \mathcal D $, the value of $ s $
increases). Assume that the equation $ Z_2 (s) = 0 $ has two roots on the
semiaxis $ s <0 $, we denote them  $ S_ {-} $ and
$ S_{+} $, $ S_ {-} \le S_ {+} $, or one root $ S_{+} $.

 In the latter case $ Z_1 (s) \to + \infty $ for $ s \to - \infty
$.

If the equation $ Z_2 (s) = 0 $  has two roots on the semiaxis
$ s <0 $, we denote them $ s_{-} $ and $ s_{+} $, $ s_{-} \le s_{ +}
$.

The trajectory $ {\mathcal D} (s, t) $ lies  between the curves $
\sqrt {Z_2 (s)} $ and $ \sqrt {Z_1 (s)} $. Denote by $ s^k_* $ the
point at which $ {\mathcal D} (s, t) $ intersects the axis $
{\mathcal D} = 0 $ for the $ k $-th time. If  $ {\mathcal D}> 0 $,
then  $ \sqrt {Z_1 (s)} \le {\mathcal D} (s, t) \le \sqrt {Z_2(s)}
$, $s$ increases. If  $ {\mathcal D}< 0 $, then $ -\sqrt {Z_1 (s)}
\le {\mathcal D} (s, t) \le -\sqrt {Z_2 (s)} $, $s$ decreases.

Let us compose the curve $L$ as follows: it consists of the part of trajectory $\sqrt{Z_2(s)}$ if
 $ {\mathcal D}> 0 $ and switches to the part of trajectory $-\sqrt{Z_1(s)}$ if  $ {\mathcal D}< 0 $
 at the axis $ {\mathcal D}= 0 $. If $-\sqrt{Z_1}(s)=0$ has two roots at $s<0$, then $L$ switches
 once again the the axis $ {\mathcal D}= 0 $ to the trajectory $\sqrt{Z_2(s)}$, and so on.
 Denote by $ L^k $ the point at which $ L $ intersects the axis $ {\mathcal D} = 0 $  for the
  $ k $-th time for $ {\mathcal D}_0<0$ and in the $ k-1 $-th time for $ {\mathcal D}_0\ge 0$.
  Thus, the trajectory  $ {\mathcal D} (s, t) $ makes at least one revolution if the number of points $ L^k $ is more
  then 2 ($L^2, L^3, \dots$) for $ {\mathcal D}_0\ge  0 $ and more then 1  ($L^1 \dots$) for $ {\mathcal D}_0< 0 $.

In its turn,  we compose the curve $l$ as follows: it consists of
the part of trajectory $\sqrt{Z_1(s)}$ if  $ {\mathcal D}> 0 $,
switches to the part of trajectory $-\sqrt{Z_2(s)}$ if  $ {\mathcal
D}< 0 $ at the axis $ {\mathcal D}= 0 $ and switches once again the
the axis $ {\mathcal D}= 0 $ to the trajectory $\sqrt{Z_1(s)}$, and
so on. Denote by $ l^k $ the point at which $ l $ intersects the
axis $ {\mathcal D} = 0 $  for the $ k $-th time for $ {\mathcal
D}_0<0$ and in the $ k-1 $-th time for $ {\mathcal D}_0\ge 0$.

 We note that $s_*^k$ for all $k \in \mathbb N$ lies between $l^k$ and $L^k$.

 \bigskip

\subsection{General estimates of $J$}


\begin{lemma}\label{L1}
\begin{equation}\label{J_gen}
- | \tilde{\bf\Xi}|^2 - \tilde {\mathcal D}^2 \le 2 J\le  {\mathcal D}^2 + |{\bf\Xi}|^2,
\end{equation}
where  $\tilde{\bf\Xi}$ and $ \tilde {\mathcal D}$ are defined as ${\bf\Xi}$ and $ {\mathcal D}$ with the change of all derivatives $\partial_i V_j$ to $|\partial_i V_j|$.
\end{lemma}
\proof The estimate \eqref{J_gen} is a corollary of the elementary inequality
$-(a-b)^2\le 2ab\le (a+b)^2, $ for every $a, b\in\mathbb R$. $\Box$

Estimate \eqref{J_gen} does not yet provide an opportunity to find explicit bounds, but we will use it for particular cases.

\subsection{Behavior in the upper hyperplane ${\mathcal D}>0$}

Now we are going to prove that if ${\mathcal D}(t)$ becomes unbounded, it happens when ${\mathcal D}<0.$
\begin{lemma}\label{L2} The projection of the trajectory ${\mathcal D}={\mathcal D}(s,t)$ to the plane $(s, {\mathcal D})$ is bounded in the upper hyperplane ${\mathcal D}>0$.
\end{lemma}

\proof As follows from the second  equation of \eqref{base4.3}, the motion along the projection ${\mathcal D}={\mathcal D}(s,t)$ onto the plane $(s, {\mathcal D })$ is clockwise for $s<0$. Therefore, in the upper hyperplane ${\mathcal D}>0$
there are two possibilities for a trajectory starting from any point $(s_0,{\mathcal D}_0>0)$:

 1. to cross the axis ${\mathcal D}=0$ and fall into the lower hyperplane, which means that the projection of the phase trajectory ${\mathcal D}(s,t)$ was bounded for ${\mathcal D }>0$,

 or

 2. to tend to $+\infty$ as $s\to -0$ and $t\to t_*\le \infty$. (This would mean that a vacuum point is formed along the characteristic).

 Assume that the latter possibility realizes, i.e. ${\mathcal D}(t)\to +\infty$ and  $s(t)\to -0$ as $t\to t_*\le \infty$.  We see from
\eqref{base4.3} that
$ \dot{\mathcal D } =   - {\mathcal D}^2 + 2J-s-1$, therefore  $ J\to +\infty$, if $t\to t_*$ and $s(t)\to -0$. Therefore $J>K>0$ for sufficiently small $s<0$. Then
\eqref{Qp} implies
\begin{eqnarray*}
 \frac{dZ}{ds} \le Q_2(s,Z)=\frac{2 (Z +s + 2 K +1)}{s},
\end{eqnarray*}
$Z_2(s) ={\mathcal D}^2 \le -2K-1-2s+C s^2$ and $Z_2(0)\le-2K-1<0$. This contradiction shows that our assumption was wrong and ${\mathcal D}>0$ is bounded. $\Box$

\subsection{Counting the number of periods for which the boundedness of density is guaranteed}\label{iter}

Based on estimates \eqref {Z_ineq1}, \eqref {Z_ineq2} we can find the number of periods for which the breaking
is guaranteed not to occur.

Let, for definiteness, $ {\mathcal D}_0 $ be negative. Assume that the initial data are such that the solution does not blow up and   on the curve $L$ there exist points $(s=L_1, {\mathcal D}=0)$ and $(s=L_2,  {\mathcal D}=0)$, $L_1<L_2$. However, we can use  $(L_2, 0)$ as initial data and check whether the equation $Z_1(s)=0$ has the second root ($L_3$) at $s<0$ besides $L_2$. This would mean that there $L$
(and  ${\mathcal D}(t,s)$) makes at least two revolutions and at the axis ${\mathcal D}=0$ there exist two points of intersection $s=L_3$ and $s=L_4$, $L_3<L_4$, where $L_3$ and $l_4$ are the roots of equation $Z_2(s)=0$. Then we can use  $(L_4, 0)$ as initial data  again and so on. We can continue this procedure until at some step $n$ for the point  $(L_{2n}, 0)$,
taken as initial data, we do not obtain a sufficient condition for boundedness of ${\mathcal D}$.
 This means that we can guarantee $n$ revolutions of the trajectory  ${\mathcal D}(t,s)$.
  The iterative process of counting revolutions is easy to implement numerically, if we know $Z_1$ and $Z_2$ (see the example in Sec.\ref{4.2}. Nevertheless,
  the estimate of numbers of revolutions from below may be rough.

\subsection{Estimates of the guaranteed time of existence of the bounded density}\label{time}

Let $ \Div \bv_0<0$ and for the characteristic starting from the
point $x_0$ we can guarantee $n$ revolutions of the trajectory
${\mathcal D}(t,s)$ around the origin. Then the guaranteed lifetime
of a solution with bounded density, denoted as $T(x_0)$, can be
estimated as
$$T_l(x_0)<T(x_0)<T_L(x_0),$$ where $T_l$ and $T_L$ is the time of passage of $n$
turns along a compound spiral lines $l$ and $L$, respectively. Thus,
\begin{eqnarray*}\label{T_l}
T_l=-\int\limits_{s_0}^{l^1}\frac{ds}{s\sqrt{Z_2(0;s)}}-\int\limits_{l^1}^{l^2}\frac{ds}{s\sqrt{Z_1(1;s)}}
-\int\limits_{l^2}^{l_3}\frac{ds}{s\sqrt{Z_2(2;s)}}-\cdots -\int\limits_{l^{2n}}^{l_{2n+1}}\frac{ds}{s\sqrt{Z_2(n;s)}},
\end{eqnarray*}
\begin{eqnarray*}\label{T_l}
T_L=-\int\limits_{s_0}^{L^1}\frac{ds}{s\sqrt{Z_1(0;s)}}-\int\limits_{L^1}^{L^2}\frac{ds}{s\sqrt{Z_2(1;s)}}
-\int\limits_{L^2}^{L_3}\frac{ds}{s\sqrt{Z_1(2;s)}}-\cdots -\int\limits_{L^{2n}}^{s_0}\frac{ds}{s\sqrt{Z_1(n;s)}},
\end{eqnarray*}
where we denoted $Z_i(k;s)$, $i=1, 2$, $k=0, 1,\dots, n$, the
functions  $Z_1(s)$ and $Z_2(s)$, based on initial
data at points $(L_{2k}, 0)$ for $Z_1(k,s)$ and $(L_{2k-1}, 0)$ for
$Z_2(k,s)$.

For $ \Div \bv_0\ge 0$ the formulas are similar, except for the sum $T_l$ begins from
 $-\int\limits_{s_0}^{l^2}\frac{ds}{s\sqrt{Z_2(0;s)}}$ and the sum $T_L$ begins from
 $-\int\limits_{s_0}^{L^2}\frac{ds}{s\sqrt{Z_1(0;s)}}$.

 The time $T_*$ of existence of the smooth solution  with a bounded density can be estimated as
 $$T_*>\inf\limits_{x_0\in \mathbb R} T_l(x_0),$$
however $\sup\limits_{x_0\in \mathbb R} T_L(x_0)$ can be less than
$T_*$.

\begin{remark}
The technique of counting the number of oscillations before the
blow-up in the case of multidimensional electrostatic
non-relativistic oscillations is similar to the case of 1D
relativistic oscillations \cite{RChZAMP21_2}.
\end{remark}

\section{Particular cases}\label{Sec3}
\subsection{1D oscillations, a criterion of a singularity formation}\label{1D}
In this case  $(\bv, \bE)=(\bv(x_1), \bE(x_1) ) $,
  $V_2=E_2=0$, ${\bf\Xi}=0$, condition \eqref{cond} evidently holds.  It is easy to see that here $J=0$
  and system
\eqref{base4.3} turns into
\begin{eqnarray*}
\dot \lambda ={\mathcal D}(1-\lambda),\quad
\dot{\mathcal D } =   - {\mathcal D}^2-\lambda.\label{base5.1}
\end{eqnarray*}
 Such a system has been considered in \cite{RChZAMP21}, where
the following criterion for the preservation of the global in time
smoothness was obtained: at each point $ x_0 \in {\mathbb R} $
(here $x=x_1$) the condition
\begin{eqnarray}\label{sufcond1D}
{ \Delta =(V_{10}')^2+2(E_{10})' - 1<0.}
\end{eqnarray}
holds.


\subsection{Globally smooth affine solutions}\label{Sec_aff_ax}

\begin{defn}
A solution $\bv, \bE$ is called  {\it an affine solution } if it has the form
\begin{equation}\label{AS}
\bv= {\mathcal Q}(t) \br,\quad \bE= {\mathcal R}(t) \br,
\end{equation}  
with $(3\times 3)$ matrices ${\mathcal Q}=(q_{ij})$ and ${\mathcal R}=(r_{ij})$.
\end{defn}
 Condition
\eqref{cond} dictates $q_{ij}=q_{ji}$,  $r_{ij}=r_{ji}$, so we get
the matrix equation
\begin{equation*}
\dot {\mathcal Q}+{\mathcal Q}^2+{\mathcal R}=0, \quad \dot {\mathcal R}=(1-{\rm tr} {\mathcal R}){\mathcal Q}
\end{equation*}
 or quadratically
nonlinear system of 12 ODEs.

For the plain oscillations $q_{3j}= q_{i3}=r_{3j}= r_{i3}=0$, $i,j=1, 2$, so the number of equations is 6.

\begin{defn}
A solution $\bv, \bE$ is called  {\it radially symmetric } if 
$$\bv=F(t,r) \br,\, \bE=G(t,r) \br,\,\br=(x_1, x_2,x_3), \, r=\sqrt{x_1^2+x_2^2+x_3^2}.$$
\end{defn}

It is easy to check that for the affine solutions and radially symmetric solutions $\Rot \bv=0$ and  condition \eqref{cond} holds.

Let us consider solutions that are both affine and radially symmetric, i.e $F=\alpha(t)$,  $G=\beta(t)$.
Here ${\mathcal D}=3\alpha$, $J=3\alpha^2$, $\lambda= 3\beta$, and the system \eqref{base4.3} takes a closed form
\begin{eqnarray*}
\dot{\alpha } =   - {\alpha}^2-\beta,\quad \dot \beta = \alpha (1-3\beta),
\label{base5.1R3}
\end{eqnarray*}
the first integral is
\begin{eqnarray*}\label{fiRad3}
\alpha^2=2\beta-1+K|1-3\beta|^\frac{2}{3} , \quad K=\frac{1-2\beta(0)+\alpha^2(0)}{|1-3\beta(0)|^\frac{2}{3}}.
\end{eqnarray*}
Recall that the positivity of density requires $ 1-3\beta>0,$ see
\eqref{n}. The prevailing term of the right hand side as $\beta \to
-\infty$ is $2\beta$, it tends to $-\infty$ for any $K$. This means
that $\alpha $ and $\beta$ are bounded {\it for any initial data},
the density is bounded and the smooth in time solution to the Cauchy
problem exists for all radially symmetric initial data.

In the case of plane oscillations we have
$$\bv=F(t,r) \br,\, \bE=G(t,r) \br,\,\br=(x_1, x_2,0), \, r=\sqrt{x_1^2+x_2^2}.$$
For the corresponding affine solution ${\mathcal
D}=2\alpha$, $J=\alpha^2$, $\lambda= 2\beta$, therefore system
\eqref{base4.3} takes a closed form
\begin{eqnarray*}
\dot{\alpha } =   - {\alpha}^2-\beta,\quad \dot \beta = \alpha (1-2\beta),
\label{base5.1R3}
\end{eqnarray*}
the first integral is
\begin{eqnarray*}\label{fiRad2}
2\alpha^2=-(1-2\beta) \ln|1-2\beta|+K(1-2\beta)  +1, \quad K=\frac{1-2 \alpha^2(0)}{1-2\beta(0)}-\ln|1-2\beta(0)|.
\end{eqnarray*}
From the positivity of density we have $ 1-2\beta>0,$ therefore the
prevailing term of the right hand side as $\beta \to -\infty$ is the
logarithmic one, it tends to $-\infty$ for any $K$. This means that
$\alpha $ and $\beta$ are bounded for any initial data, the density
is bounded and the smooth in time solution to the Cauchy problem
exists for all radially symmetric initial data (see also
\cite{CH18}, \cite{R_PhysicaD22}).

These examples show, in particular, that there are three-dimensional initial data that do not lead a finite time blow-up.
 However, as was proved in \cite{R_PhysicaD22}, radially symmetric solutions of the cold plasma equation, which are not affine, in the general case blow up even if they are arbitrarily small perturbations of the zero stationary state.

\section{The plain and irrotational oscillations }\label{Sec4}

In this section we specify two cases, where it is possible to obtain an estimate giving a possibility to study the behavior of the trajectory in the lower hyperplane.

\subsection{Plain oscillations  }\label{Sec4.1}
For the case of plain oscillations $\bv=\bv(x_1,x_2)$, $\bE=\bE(x_1,x_2)$, $V_3=E_3=0$, $\xi_1=\xi_2=0$, and we have
\begin{eqnarray*}
\dfrac{\partial \xi_3 }{\partial t}  &+& ( \bv \cdot \nabla)  \xi_3  = - {\mathcal D}\xi_3.
\end{eqnarray*}
Along characteristics this equations can be written as $\dot { \xi_3} =- {\mathcal D}{ \xi_3}$,
therefore  \eqref{base4.3} implies
\begin{eqnarray}\label{FE}
\xi_3= C_3 (\lambda-1),\quad C_3={\rm const}.
\end{eqnarray}
Thus, from \eqref{J_gen} we have
\begin{eqnarray*}
 2 J\le {\mathcal D} ^2 + C^2_3 s^2, 
\end{eqnarray*}
and \eqref{Qm} takes the form
\begin{eqnarray*}\label{Qm2D}
 Q_1(s,Z) = \frac{2 (2 Z +s + C_3^2 s^2 +1)}{s}.
\end{eqnarray*}
Therefore
\begin{eqnarray*}
Z_1(s)=A_4 s^4- A_2 s^2 +A_1 s+A_0,
\end{eqnarray*}
where
\begin{eqnarray*}\label{ZpA1}
A_4= \frac{Z_0+{\xi}_{30}^2+\frac{2}{3}s_0 +\frac12}{s_0^4}, \quad A_2=-\frac{|{\xi}_{30}|^2}{s_0^2},\quad  A_1=-\frac23, \quad A_0=-\frac12,\label{ZpA2}
\end{eqnarray*}
the coefficients are substituted with the value of the constant $ C_3 $, found from \eqref{FE}.



\bigskip

\subsection{Irrotational oscillations }\label{Sec4.2}

The case of irrotational oscillations is analogous.
Since
\begin{eqnarray*}
\dfrac{\partial \xi_i }{\partial t}  &+& ( \bv \cdot \nabla)  \xi_i  = - {\mathcal D}\xi_i-({\bf \Xi} \cdot \nabla ) V_i,\quad i=1, 2, 3,
\end{eqnarray*}
then the uniqueness of the solution of the Cauchy problem implies that if ${\bf\Xi}_0=0$, then ${\bf\Xi}=0$ identically. Therefore,
 \eqref{J_est} gives
\begin{eqnarray*}
 2 J\le {\mathcal D} ^2 \label{J1}
\end{eqnarray*}
and \eqref{Qm} is
\begin{eqnarray}\label{QmIrr}
 Q_1(s,Z) = \frac{2 (2 Z +s+1)}{s}.
\end{eqnarray}
Thus,
\begin{eqnarray*}
Z_1(s)=A_4 s^4+A_1 s+A_0,
\end{eqnarray*}
where
\begin{eqnarray*}
A_4= \frac{Z_0+\frac{2}{3}s_0 +\frac12}{s_0^4}, \quad  A_1=-\frac23, \quad A_0=-\frac12.
\end{eqnarray*}


For the both case the following lemma holds.

\begin{lemma} \label{L3}
Let at a point $ x_0 \in {\mathbb R} ^ 3 $ the initial conditions
\eqref{CD1} are such that $\bv_0=(V_1(x_1, x_2), V_2(x_1, x_2), 0),$ $\bE_0=(E_1(x_1, x_2), E_2(x_1, x_2), 0),$   or  ${\rm curl} \bv_0=0$ and
\begin{eqnarray}\label{sufcond0}
{\rm div} \bv_0<0, \quad { \Delta_- =({\rm div} \bv_0)^2+|{\rm curl} \bv_0|^2+\frac{2}{3}{\rm div} \bE_0 -\frac16< 0,}
\end{eqnarray}
or
\begin{eqnarray}\label{sufcond01}
{\rm div} \bv_0=0, \quad {\rm div} \bE_0>0,\quad { \Delta_-=|{\rm curl} \bv_0|^2+\frac{2}{3}{\rm div} \bE_0 -\frac16 < 0}.
\end{eqnarray}
Then there exists a moment $t_*$ such that the  trajectory $(s, \mathcal D)$ turns out in the upper hyperplane ${\mathcal D}<0$
(a breaking of oscillations does not occur in
the first period).
\end{lemma}

\proof The sense of Lemma \ref{L1} is such that the  trajectory ${\mathcal D}(s,t)$ does not go to infinity during the first oscillation
in the lower hyperplane ${\mathcal D}<0$.
 As follows from Sec.\ref{Sec2}, the trajectory in the lower hyperplane ${\mathcal D}<0$ is bounded from below by
 the curve $-\sqrt{Z_2(s)}$. In its turn, ${Z_2(s)}$ is bounded if and only if $A_4<0$, where $A_4$ is defined in \eqref{ZpA1}.
 The condition $A_4<0$ is exactly \eqref{sufcond0}. Further, since the motion along the phase trajectory is clockwise, then
 \eqref{sufcond01} signifies that the trajectory will occur in the lower hyperplane with the condition \eqref{sufcond0} at moment $t>0$.
  $\Box$

Thus, as a corollary of Lemmas \ref{L2} and \ref{L3}, we get the following theorem.

 \begin{theorem} \label{T1}
Let at each point $ x_0 \in {\mathbb R} ^ 3 $ the initial conditions
\eqref{CD1} are such that $\bv_0=(V_1(x_1, x_2), V_2(x_1, x_2), 0),$ $\bE_0=(E_1(x_1, x_2), E_2(x_1, x_2), 0),$  or  ${\rm  rot} \bv_0=0$ and one of conditions
\eqref{sufcond0} and \eqref{sufcond01} holds.
Then the density $n$,  found as \eqref{n} from
the solution of the problem \eqref {4}, \eqref{CD1}, is bounded 
during the first period of oscillations.
\end{theorem}

\proof According to conditions \eqref{sufcond0} or \eqref{sufcond01}, every  trajectory ${\mathcal D}={\mathcal D}(s,t)$ starts from a point do not laying in the upper hyperplane and turns out in the upper hyperplane (Lemma \ref{L3}) and then again turns out in the lower hyperplane (possibly, at the point, do not satisfying to  \eqref{sufcond0} or \eqref{sufcond01}) (Lemma \ref{L2}). However, the trajectories starting from each point $ x_0 \in {\mathbb R} ^ 3 $ make as least one full rotation around the origin $\lambda={\mathcal D}=0$, therefore ${\mathcal D}$ is bounded during the first period of oscillations, and \eqref{n} implies that $n$ is also bounded. $\Box$


\subsection{Affine solutions, an example}\label{Sec4.3}

 For the affine solutions \eqref{AS} Theorem \ref{T1} gives
the following sufficient conditions for maintaining smoothness in the
first period of oscillations:
\begin{equation*}\label{sufcond_aff}
{\rm tr}\, Q(0)<0,\quad \Delta_- =({\rm tr} \, Q(0))^2+\frac23 {\rm tr} \, R(0)-\frac16< 0,
\end{equation*}
or
\begin{equation}\label{sufcond_aff2}
{\rm tr}\, Q(0)= 0,\quad {\rm tr}\, R(0)>0, \quad \Delta_- =4 \,{\rm tr}\, R(0)-1< 0.
\end{equation}
For the axisymmetric 3D case,  \eqref{sufcond_aff2} implies
 $\beta< \frac{1}{12}$, whereas as is shown in Sec.\ref{Sec_aff_ax} the solution is globally smooth for $\beta< \frac13$.
 
 For the axisymmetric 2D and 1D cases,  \eqref{sufcond_aff2} gives
 $\beta< \frac18$ and $\beta< \frac14$, respectively. However,  the solution is globally smooth for $\beta<\frac12$ for both these cases.

\section{Radially symmetric oscillations in any dimensions}\label{Sec5}

For the radially symmetric case we can obtain a two-sided estimate \eqref{J_est} and obtain the find the time for which we guarantee the existence of the smooth solution to the problem \eqref{4}, \eqref{CD1}. This result can be considered as an addition to the results of the recent paper \cite{R_PhysicaD22}, where we proved that any solution to the Cauchy problem for $\bd=2$ and $\bd=3$ generally blows up, even for arbitrary small perturbation to the zero steady state.

The calculations are similar in the space of any dimension $\bd\ge 1$,
we assume
$$\bv=F(t,r) \br,\, \bE=G(t,r) \br,\,\br=(x_1, x_2, \dots, x_\bd), \, r=\sqrt{x_1^2+x_2^2+\dots +x_\bd^2}.$$

However, we are interested mostly  in physical dimension $\bd=3$. For the plain oscillations ($\bd=2$ )
\begin{equation}\label{sol_form}
  \bv=F(t,r) \br,\, \bE=G(t,r) \br,\,\br=(x_1, x_2, 0), \, r=\sqrt{x_1^2+x_2^2}.
\end{equation}
It is easy to check that $\Rot \bv=0$ and  condition \eqref{cond} holds.

Further,
\begin{eqnarray*}
 &&{\mathcal D} = \bd \, F +F_r \, r, \quad  \lambda = \bd\, G +G_r \,
 r,
 \quad J=(\bd-1)\,F \,F_r\,r+\frac{(\bd-1)\bd}{2}\, F^2,
\end{eqnarray*}
therefore
\begin{eqnarray*}
  2J &=&2\,(\bd-1)\,{\mathcal D}\, F-(\bd-1)\,\bd \,F^2.
\end{eqnarray*}

Further, if we substitute \eqref{sol_form} to \eqref{4}, we get
 \begin{eqnarray}\label{GF}
  \dot{G}
   =F-\bd F G,\quad  \dot{F}
   =-F^2 - G,
 \end{eqnarray}
 where $\dot f= \dfrac{\partial f}{\partial t}+F r \dfrac{\partial f}{\partial r}$.

 On the phase plane $(G,F)$ system \eqref{GF} implies one equation
\begin{eqnarray*}
\frac12 \dfrac{d F^2}{d G}=-\frac{ F^2 +
    G}{1-\bd
    G},
\end{eqnarray*}
which is linear with respect to $Y=F^2$ and can be explicitly integrated. Indeed,
we have for $\bd=2$
\begin{eqnarray}\label{fiRad2}
2 Y=(2G-1) \ln|1-2G|+C_2(2G-1)  -1, \\ C_2=\frac{1+2
F^2(0,r_0)}{2G(0,r_0)-1}-\ln|1-2G(0,r_0)|,\nonumber
\end{eqnarray}
for $\bd=1$ and $\bd\ge 3$
\begin{eqnarray}\label{fiRad3}
&&Y=\frac{2 G-1}{\bd-2}+C_\bd |1-\bd G|^\frac{2}{\bd} , \quad
C_\bd=\frac{1-2 G(0,r_0)+(\bd-2) F^2(0,r_0)}{(\bd-2)|1-\bd
G(0,r_0)|^\frac{2}{\bd}}.
\end{eqnarray}

 If $G(0,r_0)<\frac{1}{\bd}$, $r_0\in \bar {\mathbb R}$, then
the phase trajectories of \eqref{GF} are in the half-plane
$G<\frac{1}{\bd}$. In this half-plane the leading
term in the right hand side of \eqref{fiRad2} is $(2\,G-1)
\ln|1-2\,G|$, and the leading term in the right hand side of
\eqref{fiRad3} is $\frac{2\, G-1}{\bd-2}$.  Therefore $F$ and $G$
are bounded for any $t>0$.

Along every characteristic curve starting from $r_0\in\bar{\mathbb R}_+$ the functions
 $F(t)$ and $G(t)$, solution of \eqref{GF}, are periodic with the period
\begin{equation*}\label{T}
  T= 2 \int\limits_{G_-}^{G_+} \frac {d\eta}{(1-\bd\, \eta )F(\eta)},
\end{equation*}
$F$ is given as \eqref{fiRad2} or  \eqref{fiRad3}, $G_-<0$ and
$G_+>0$ are the lesser and greater roots of the equation $F(G)=0$.
Moreover, $\int\limits_0^{T} F(\tau) \, d\tau=0$.

We refer to \cite{R_PhysicaD22} for a detailed analysis of the phase portrait and the period of oscillations.

\bigskip
Thud, we can obtain the evaluating functions
$Q_1$ and $Q_2$ in \eqref{eqtQ}.

 Indeed,
  the following upper estimate holds:
 \begin{eqnarray}\label{est_Jp}
   2J\le 2 (\bd-1) F \mathcal D \le\\ -\left(\sigma_1 \mathcal D-\frac{(\bd-1)F}{\sigma_1}\right)^2 +\sigma_1^2 \mathcal D^2+\frac{(\bd-1)^2 }{\sigma_1^2}\, F^2\le \sigma_1^2 \mathcal D^2+\frac{(\bd-1)^2 }{\sigma_1^2}\, F_+^2,\quad \sigma_1>0.\nonumber
 \end{eqnarray}
The lower estimate is analogous:
\begin{eqnarray}\label{est_Jm}
   2J\ge 2(\bd-1) F \mathcal D-  \bd (\bd-1) F_+^2 \ge\\
   \left(\sigma_2 \mathcal D+\frac{(\bd-1) F}{\sigma_2}\right)^2 -\sigma_2^2 \mathcal D^2-\frac{(\bd-1)^2 F^2}{\sigma_2^2}- \bd (\bd-1) F_+^2\ge\nonumber\\
   -\sigma_2^2 \mathcal D^2-\frac{(\bd-1)(\bd (\sigma_2^2+1)-1)}{\sigma_2^2}\,F_+^2,\quad
   \sigma_2>0.\nonumber
 \end{eqnarray}
 The parameters $\sigma_i$, $i=1,2$, in both estimates can be different and are chosen at our convenience.

 Thus, taking into account \eqref{est_Jp} and \eqref{est_Jm}, we get
 estimating functions
 \begin{eqnarray*}\label{Qm1}
 \bar Q_1(s,Z) = \frac{2 ((1+\sigma_1^2) Z +s +\frac{(\bd-1)^2 F_+^2}{\sigma_1^2} +1)}{s},
\end{eqnarray*}
\begin{eqnarray*}\label{Qp1}
\bar Q_2(s,Z) = \frac{2  ((1-\sigma_2^2) Z +s -\frac{(\bd-1)(\bd (\sigma_2^2+1)-1) F_+^2}{\sigma_2^2} +1)}{s}.
\end{eqnarray*}

  The equations $ \frac{d Z}{d s}=\bar Q_1 (s,Z)$ and $ \frac{d Z}{d s}=\bar Q_2 (s,Z)$ can be solved, their solutions are
\begin{eqnarray}\label{Zp1}
&&\bar Z_1(s)= -\frac{2s}{1+2\sigma_1^2}-\frac{(\bd-1)^2 F_+^2+\sigma_1^2}{\sigma_1^2(1+\sigma_1^2)}+C_1\,s^{2(1+\sigma_1^2)},\\
&&C_1= s_0^{-2(1+\sigma_1^2)}\left(Z(0)+ \frac{2
\sigma_1^2(1+\sigma_1^2) s_0+(1+2 \sigma_1^2)
((\bd-1)^2 F_+^2+\sigma_1^2)}{\sigma_1^2
(1+2\sigma_1^2)(1+\sigma_1^2)}\right),\label{C1}
\end{eqnarray}
and
\begin{eqnarray}\label{Zm1}
&&\bar Z_2 (s)= -\frac{2s}{1-2\sigma_2^2}+\frac{K-1}
{1-\sigma_2^2}+C_2\,s^{2(1-\sigma_2^2)},\quad K=\frac{(\bd-1)(\bd (\sigma_2^2+1)-1) F_+^2}{\sigma_2^2},
\\&&C_2= s_0^{-2(1-\sigma_2^2)}\left(Z(0)+\frac{2 (s_0-1)}{2 \sigma_2^2-1}+\frac{K}{1-\sigma_2^2}
\right).\label{C2}
\end{eqnarray}

Thus, we get
\begin{equation*}
   \bar Q_1(s,Z) \le Q(s,Z,J)\le  \bar Q_2(s,Z),
\end{equation*}

The new estimate $
   \bar Q_1(s,Z) \le Q(s,Z,J) $
is not necessarily better that the previous estimate  $
    Q_1(s,Z) \le Q(s,Z,J) $, see \eqref{QmIrr}.
Nevertheless, for the initial data \eqref {CD1} sufficiently small
in the uniform norm, the new estimate improves the lower bound for
the existence time of a smooth solution.

\section{A  model example in 2D}\label{Sec6}

We restrict ourselves to a special kind of
initial data, the most interesting from the point of view of
physics.

Motivated by the form of a standard laser pulse \cite{Shep13} we
choose as the initial data \eqref{CD1}
\begin{equation*}\label{gauss}
|\bE_0(\rho)| = \left(\dfrac{a_*}{\rho_*}\right)^2 \rho \exp\left(-
\dfrac{\rho^2}{\rho_*^2}\right),\quad \bV_{0}(\rho) = 0, \quad
\rho=\sqrt{x_1^2+x_2^2},
\end{equation*}
$a_*$ and $\rho_*$ are parameters. For the sake of simplicity, we
change the space variable to $r=\rho/\rho_*$ and reduce the data to
\begin{equation}\label{gauss}
|\bE_0(r)| = K  \,e^{-r^2}\,{\bf r},\quad \bV_{0}(r) = 0, \quad
r=\sqrt{x_1^2+x_2^2}, \quad K=\dfrac{a_*^2}{\rho_*}>0.
\end{equation}

\begin{center}
\begin{figure}[htb]
\hspace{-1.5cm}
\begin{minipage}{0.4\columnwidth}
\includegraphics[scale=0.4]{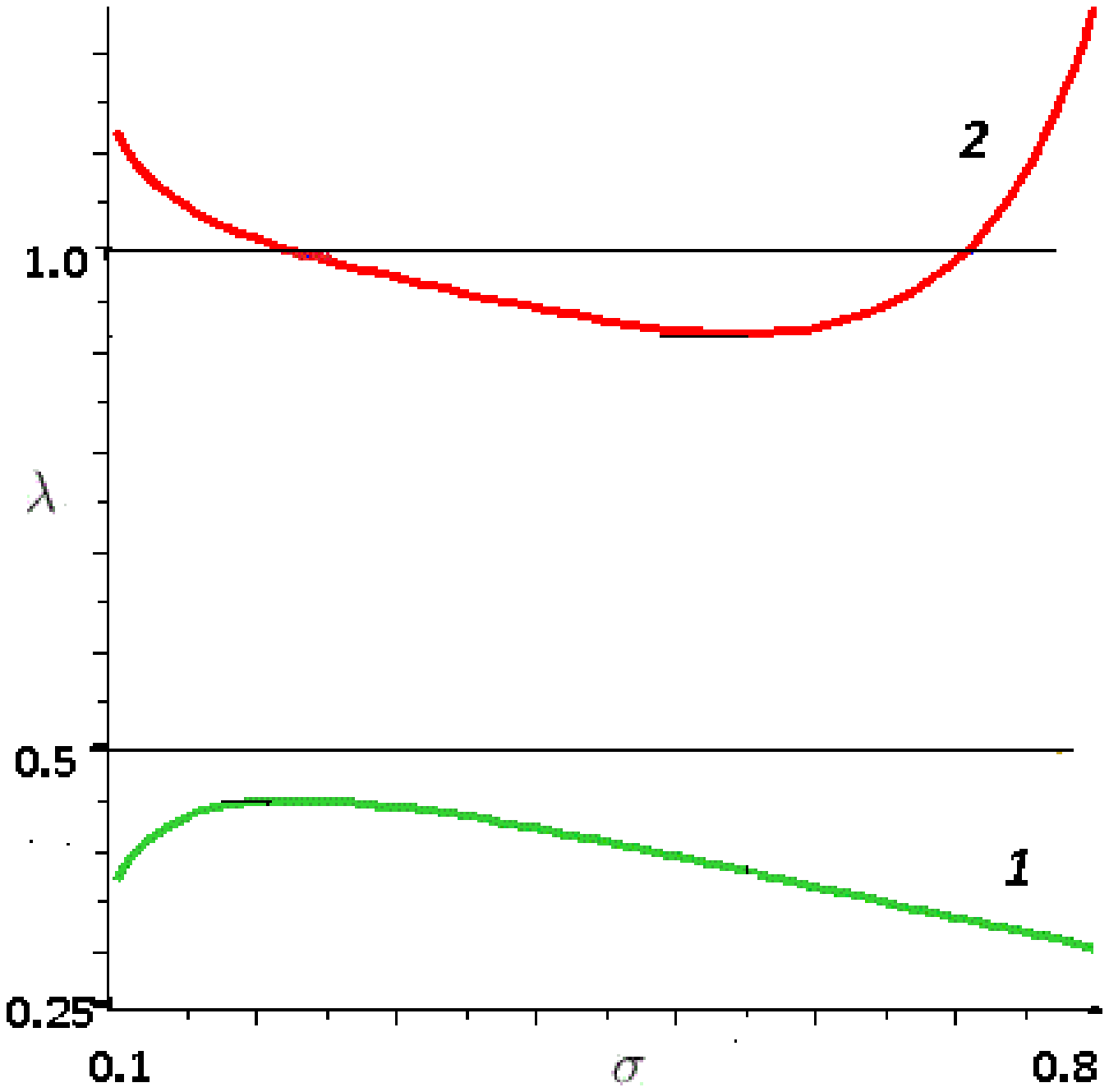}
\end{minipage}
\hspace{1.5cm}
\begin{minipage}{0.4\columnwidth} 
\includegraphics[scale=0.4]{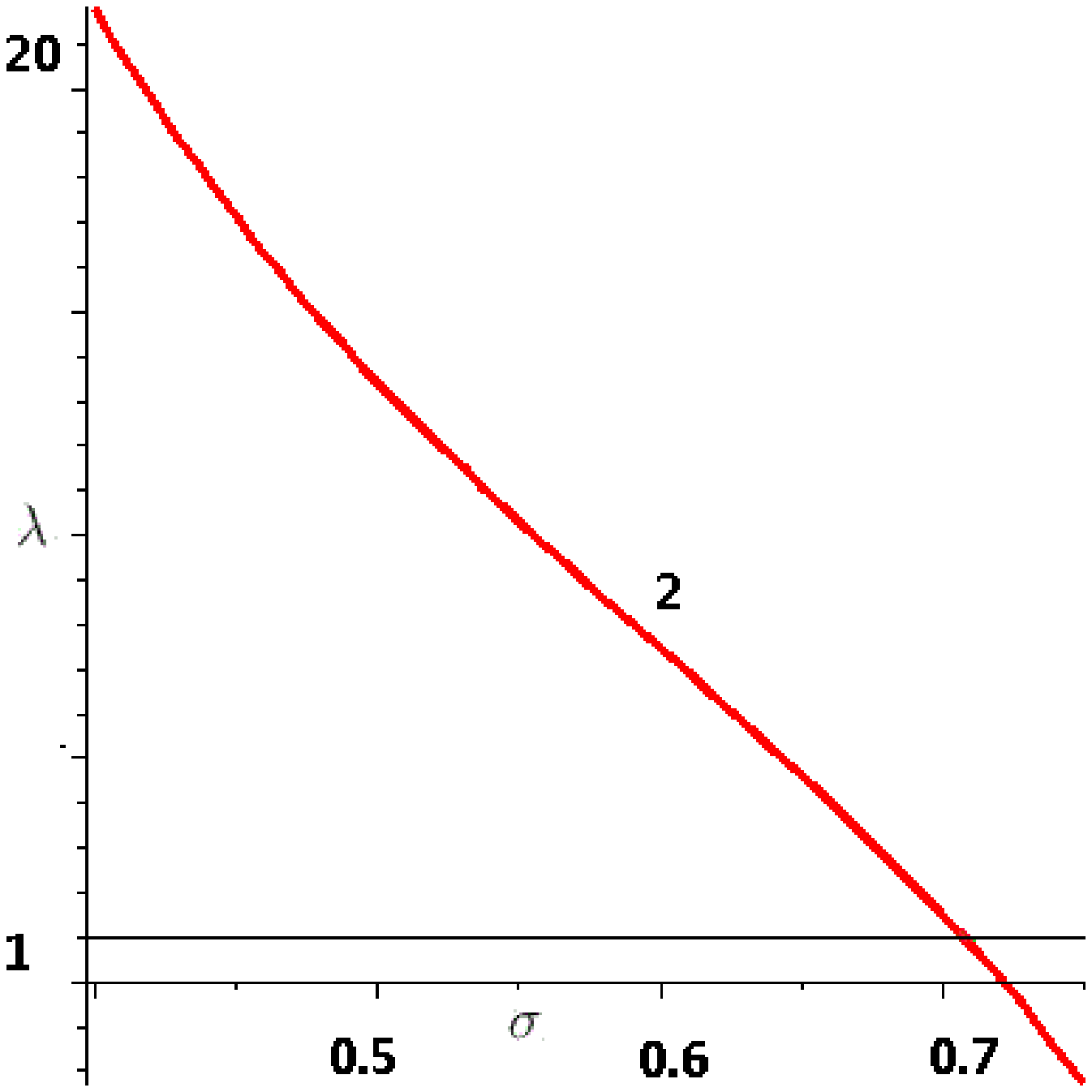}
\end{minipage}
\caption{Graphs of $C_1(\lambda, \sigma)=0$ (1) and
 $C_2(\lambda, \sigma)=0$ (2), $\lambda=s+1$, for $F_+=0.05$,
 $F_+<\frac{\sigma^4}{2\sigma^2+1}$, left.
Graph of
 $C_2(\lambda, \sigma)=0$ (2),  for $F_+=2.5$,
 $F_+>\frac{\sigma^4}{2\sigma^2+1}$, right.
  }\label{FPic}
\end{figure}
\end{center}

Fig.\ref{FPic} presents the curves $C_1=0$ and $C_2=0$ on the plane
$(\sigma, \lambda)$  ($\sigma=\sigma_1 $ or $\sigma_2$) for $\mathcal D_0=0$ for different values of
$F_+$. Below (above) $C_i$, lie such values of $\lambda_0=s_0+1$
that $C_i<0$ ($C_i\ge0$) and $Z_i$ is bounded (unbounded), $i=1,2$.
This conclusion follows from the analysis of leading terms in
\eqref{Zp1}, \eqref{Zm1}.

Since the
trajectory on the phase plane can go to infinity only at $\mathcal
D<0$, we study the estimate function under this condition.

\subsection{$\mathcal D<0$, above estimate $\bar Z_1$}
\bigskip
Thus, we have to study the curve $C_1=0$, given as \eqref{C1} for
$Z(0)=0$.  We find $s_0$ from this equation and  denote
\begin{equation*}
   S_1=-\frac12
  \frac{(1+2\sigma_1^2)(F_+^2+\sigma_1^2)}{\sigma_1^2 (1+\sigma_1^2)}.
  \end{equation*}
  Condition \eqref{sufcond01} implies that that the solution keeps
  smoothness at the first rotation provided $\lambda_0<\frac14$ ($s_0<-\frac34$) at the most "dangerous" point $r=0$
  ($G(r)$ has the maximum here).

\begin{figure}[htb]
\hspace{-1.5cm}
\begin{minipage}{0.4\columnwidth}
\includegraphics[scale=0.4]{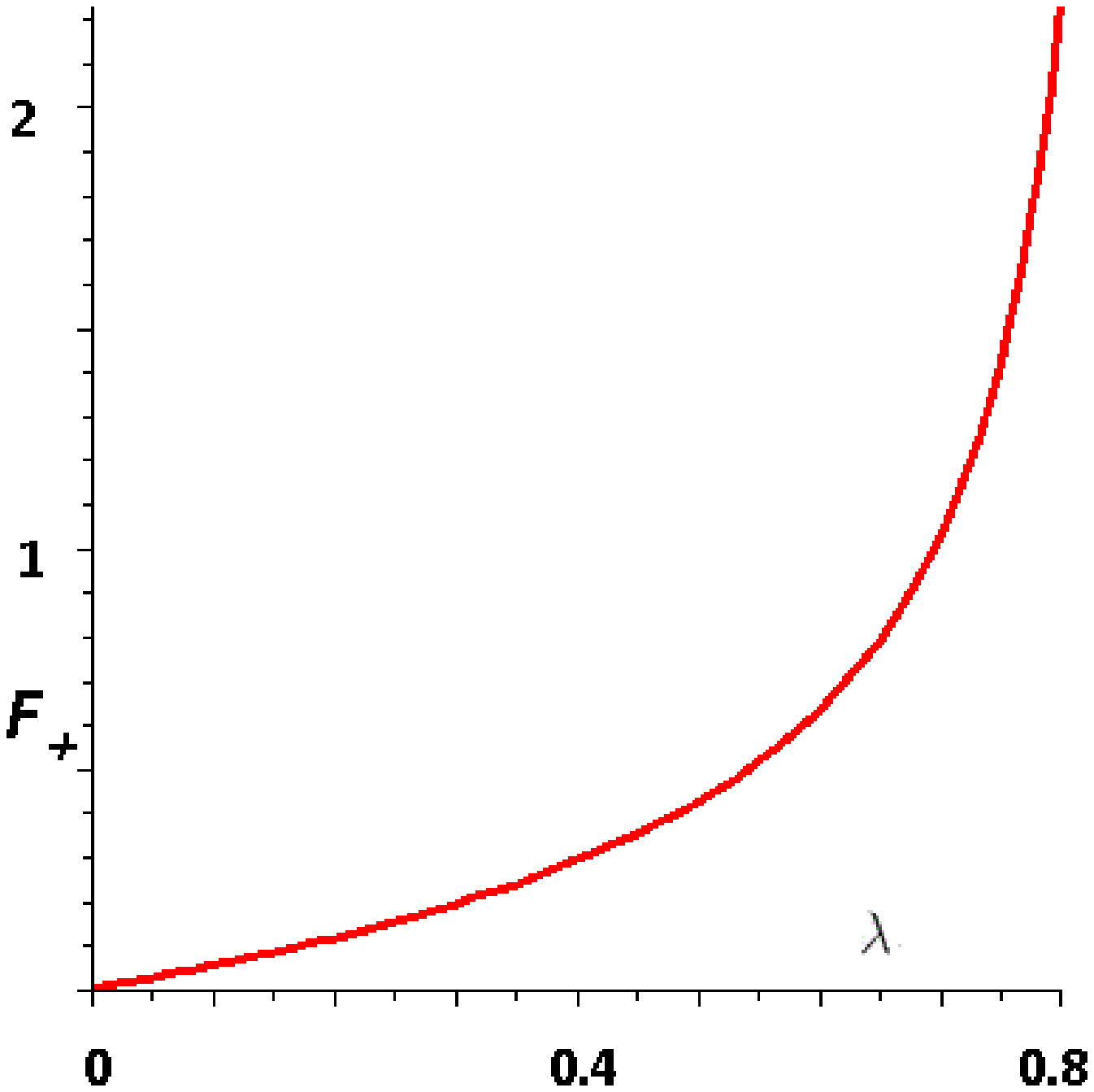}
\end{minipage}
\hspace{1.5cm}
\begin{minipage}{0.4\columnwidth} 
\includegraphics[scale=0.4]{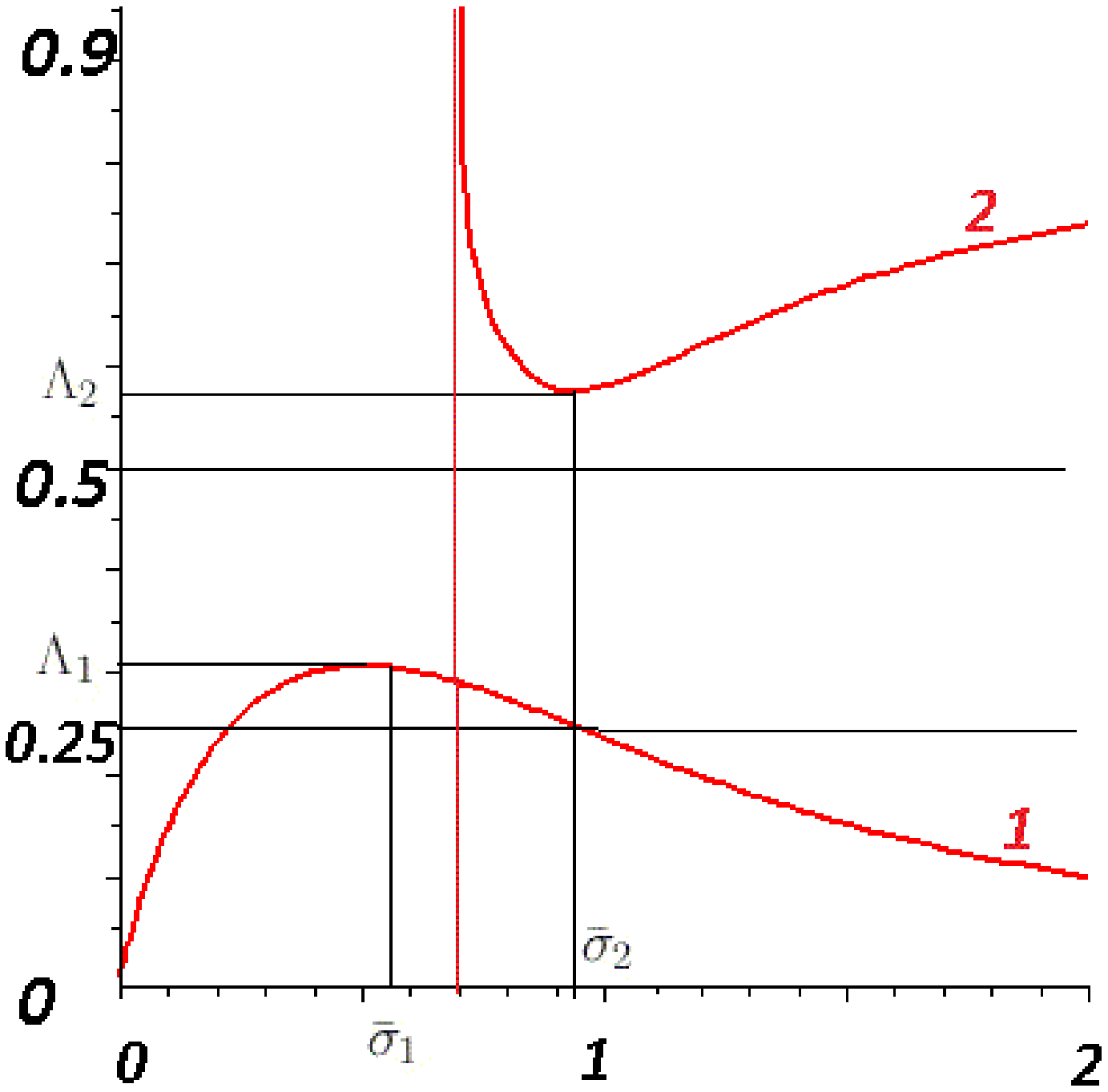}
\end{minipage}
\caption{Dependency $F_+(\lambda)$ (see \eqref{Fp}), left.
Dependency $\lambda^*_1(\sigma)$ (1) and $\lambda^*_2(\sigma)$ (2);
the maximum $(\bar \sigma_1, \Lambda_1)$ and the minimum $(\bar
\sigma_2, \Lambda_2)$, right. }\label{FS}
\end{figure}

 We have to keep in mind that $F_+$ cannot be chosen arbitrary.
Indeed, according to \eqref{sol_form}, at the point $r=0$ we chose
$\lambda_0$ together with $G(0)=\frac12 \lambda_0$.

Then we can use \eqref{fiRad2}, where the constant $C_2$ is found from
$Y(0)=0$, $G(0)=\frac12 \lambda_0$.
 Then
 \begin{equation}\label{Fp}
 \bar F_+= F_+(\lambda_0)=\frac12
\sqrt{4e^{-C-1}-2},\quad C=\frac{1}{\lambda_0-1} -\ln \left(
\frac{1-\lambda_0}{2}\right),
 \end{equation}
where $F_+(\lambda_0)=\sqrt{Y(G_m)}$, $ G_m=-e^{-C-1}$ is the
maximum point  of $Y(G)$. Graph of $F_+(\lambda_0)$ is presented in
Fig.\ref{FS}, left.

Let us introduce the function $\lambda_1(F_+(\lambda_0),\sigma_1)=
S_1(F_+(\lambda_0),\sigma_1)+1 \, \lambda_0\in (0,1)$. Taking into
account \eqref{Fp}, we can write an explicit expression for this
function
\begin{equation*}
\lambda_1(F_+(\lambda_0),\sigma_1)=
\frac{1+4\sigma_1^2-(2\sigma_1^2+1)(1-\lambda_0)
e^\frac{\lambda_0}{1-\lambda_0} }{4\sigma_1^2(\sigma_1^2+1)}.
\end{equation*}

We look for $\lambda_0\le\lambda_1(\lambda_0)$. If such
$\lambda_0>\frac14$, we obtain a sufficient condition for the
smoothness of the solution  on the first oscillation. Thus, we have
to find the fixed point $\lambda_{1}^*$ of the mapping
$\lambda_1(F_+(\lambda_0))$. This value can be expressed in terms of
the Lambert W function \cite{Lambert}:
\begin{eqnarray*}
\lambda^{*}_1(\sigma_1)&=&\frac{(4\sigma_1^4-1) L_W(\sigma_1)
 -4\sigma_1^2-1}{(4\sigma_1^4-1) \,
 L_W(\sigma_1)-4\sigma_1^2(1+\sigma_1^2)}
,\\
L_W&=&{\rm LambertW}\left(k,\frac{1}{2\sigma_1^2-1}\,
 e^{\frac{4\sigma_1^2+1}{4\sigma_1^2-1}}
 \right),
\end{eqnarray*}
where $k=-1$ for $\sigma_1<\frac{1}{\sqrt(2)}$  and $k=0$ for
$\sigma_1>\frac{1}{\sqrt(2)}$.

 Let us maximize
$\lambda^{*}_1(\sigma_1)$ (this can be done numerically). The
maximum point is $\bar\sigma_1=0.5032\dots$, the maximum value
$\Lambda_{1}\equiv\lambda^{*}_1(\bar\sigma_1)=0.3058\dots$. We can
see that the sufficient condition \eqref{sufcond01} is improved,
because now we can guarantee at least one revolution of every
trajectory on the phase plane for the data \eqref{gauss} with
$\lambda_0\in [0, \lambda_*^1]$, i.e. $K< \frac{\lambda_*^1}{2}=
0.1529\dots$ instead of $K< 0.125$.

 \subsection{$\mathcal D<0$, below estimate $\bar Z_2$}

 The analysis  is analogous to the previous
 subsection.  We find $s_0$ from  the equation $C_2=0$,  given as \eqref{C2} for $Z(0)=0$, and
   denote
   \begin{equation}\label{S2}
   S_2=\frac12
\frac{(2\sigma_2^2-1)(F_+^2 (2\sigma_2^2+1)-\sigma_2^4)}{\sigma_2^2
(\sigma_2^2-1)}.
\end{equation}
Since the denominator of $S_2$ vanishes at $\sigma_2=1$, we restrict
ourselves by the interval $\sigma_2\in (0,1)$.

For $F_+<\bar F_+=\frac{\sigma_2^4}{2\sigma_2^2}$ the function $S_2$
has a minimum of $S_2$ on the interval $\sigma_2\in (0,1)$, for
$F_+>\bar F_+$ the function $S_2$ decays with $\sigma_2$ (see
Fig.\ref{FPic}).

We consider the function $\lambda_2(F_+(\lambda_0),\sigma_2)=
S_2(F_+(\lambda_0),\sigma_2)+1 \, \lambda_0\in (0,1)$, its explicit
expression is
\begin{equation*}
\lambda_2(F_+(\lambda_0),\sigma_2)=\frac{(2\sigma_2^2-1)(\frac12
((1-\lambda_0) e^\frac{\lambda_0}{1-\lambda_0}-1)
(2\sigma_2^2+1)-\sigma_2^4)}{2\sigma_2^2(\sigma_2^2-1)}.
\end{equation*}

Now we need to find $\lambda_0\ge\lambda_2(\lambda_0)$. If such
$\lambda_0<1$, we obtain a sufficient condition for the blow-up on
the first oscillation. The fixed point $\lambda_{2}^*$ of the
mapping $\lambda_2(F_+(\lambda_0))$ exists only for
$\sigma_2>\frac{1}{\sqrt{2}}$, it can be expressed in terms of the
Lambert W function:
\begin{eqnarray*}
\lambda^{*}_2(\sigma_2)&=&\frac{(2\sigma_2^2-1) \Sigma_1(\sigma_2)
L_W(\sigma_2)+\Sigma_2(\sigma_2)}{(2\sigma_2^2-1) \Sigma_1(\sigma_2)
L_W(\sigma_2)-4\sigma_2^2(\sigma_2^2-1)},\\
 \Sigma_1&=&2\sigma_4+2\sigma_2+1,\quad
 \Sigma_1=4\sigma_6-2\sigma_4+4\sigma_2-1,\\
 L_W&=&{\rm LambertW}\left(-1,\frac{1+2\sigma_2^2}{\Sigma_1}\,
 e^{\frac{\Sigma_2}{(1-2\sigma_2^2)\Sigma_1}}
 \right).
\end{eqnarray*}
We minimize $\lambda^{*}_2(\sigma_2)$ with respect to $\sigma_2$.
The minimum point is $\bar\sigma_2=0.9423\dots$, the minimum value
$\Lambda_{2}\equiv\lambda^{*}_2(\bar\sigma_2)=0.5754\dots$.

Fig.\ref{FS}, right, presents the functions
$\lambda^{*}_1(\sigma_1)$, $\lambda^{*}_2(\sigma_2)$ and their
extrema.

We summarize our results.

\begin{prop}\label{prop}
Let us consider the solution to the Cauchy problem \eqref{4},
\eqref{gauss}.
\begin{itemize}
\item If $K< \frac12 \Lambda_1=0.1529\dots $, then the solution keeps
$C^1$ - smoothness during at least the first oscillation;
\item
if $K> \frac12 \Lambda_2=0.2877\dots $, then the solution blows up
within the first oscillation.
\end{itemize}
\end{prop}

\bigskip


\subsection{Example of the estimate of the number of oscillations before the blow up}\label{4.2}

 Let us choose the data \eqref{gauss} with $K=0.1$.
 We use the technique, described in
Secs. \ref{iter} and \ref{time} by means
of functions $\bar Z_1$ and $\bar Z_2$.

Let us list its steps.

1. Given $K=\frac{1}{2}\lambda_0$, $K<\frac{1}{2}\Lambda_1$,
 find $F_+(\lambda_0)$ by formula \eqref{Fp}.

2. For this value of $F_+$ we consider $\bar Z_1=\bar Z_1(\lambda,
\bar \sigma_1)$ (see \eqref{Zp1}).  The graph of $-\sqrt{\bar Z_1}$
is the first part of the curve $L$.

3. We find $\lambda_{0}^1<0$, the second root of equation $\bar
Z_1(\lambda,\bar\sigma_1)=0$, and use it as the new initial data.
Thus, we find new $F_+$ and  construct the function $\bar
Z_2(\lambda,\bar\sigma_2)$.  The graph of $\sqrt{\bar Z_2}$ is the
next link of the curve $L$.

4. We find $\lambda_{0}^2>0$, the second root of equation $\bar
Z_2(\lambda,\bar\sigma_2)=0$, and use it as the new initial data,
and repeat Step 2. Thus, we get the third link of the curve $L$, and
so on.

5. To find the guaranteed number of revolutions ($n$) we construct
the curve $L$ as described in Sec.\ref{iter}. The iteration
procedure stops if on the next step $C_1\ge 0$.

6. We construct the curve $l$ as described in Sec.\ref{iter} for
this $n$. The constant $T_l$, computed by the formula from
Sec.\ref{time}, is the estimate of the time of the existence of the
smooth solution from below.

\begin{figure}[htb]
\hspace{-1.5cm}
\begin{minipage}{0.4\columnwidth}
\includegraphics[scale=0.4]{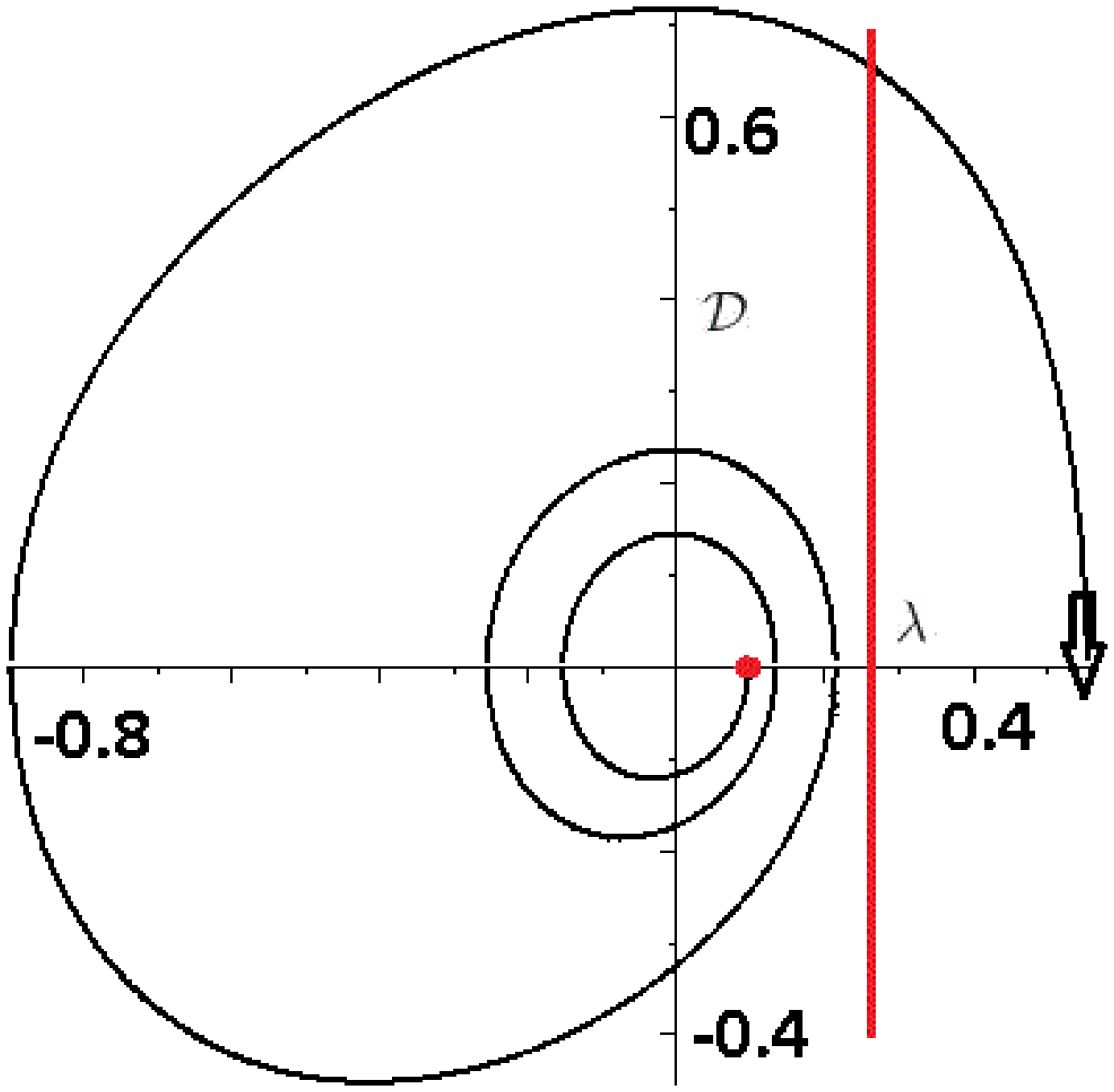}
\end{minipage}
\hspace{1.5cm}
\begin{minipage}{0.4\columnwidth} 
\includegraphics[scale=0.4]{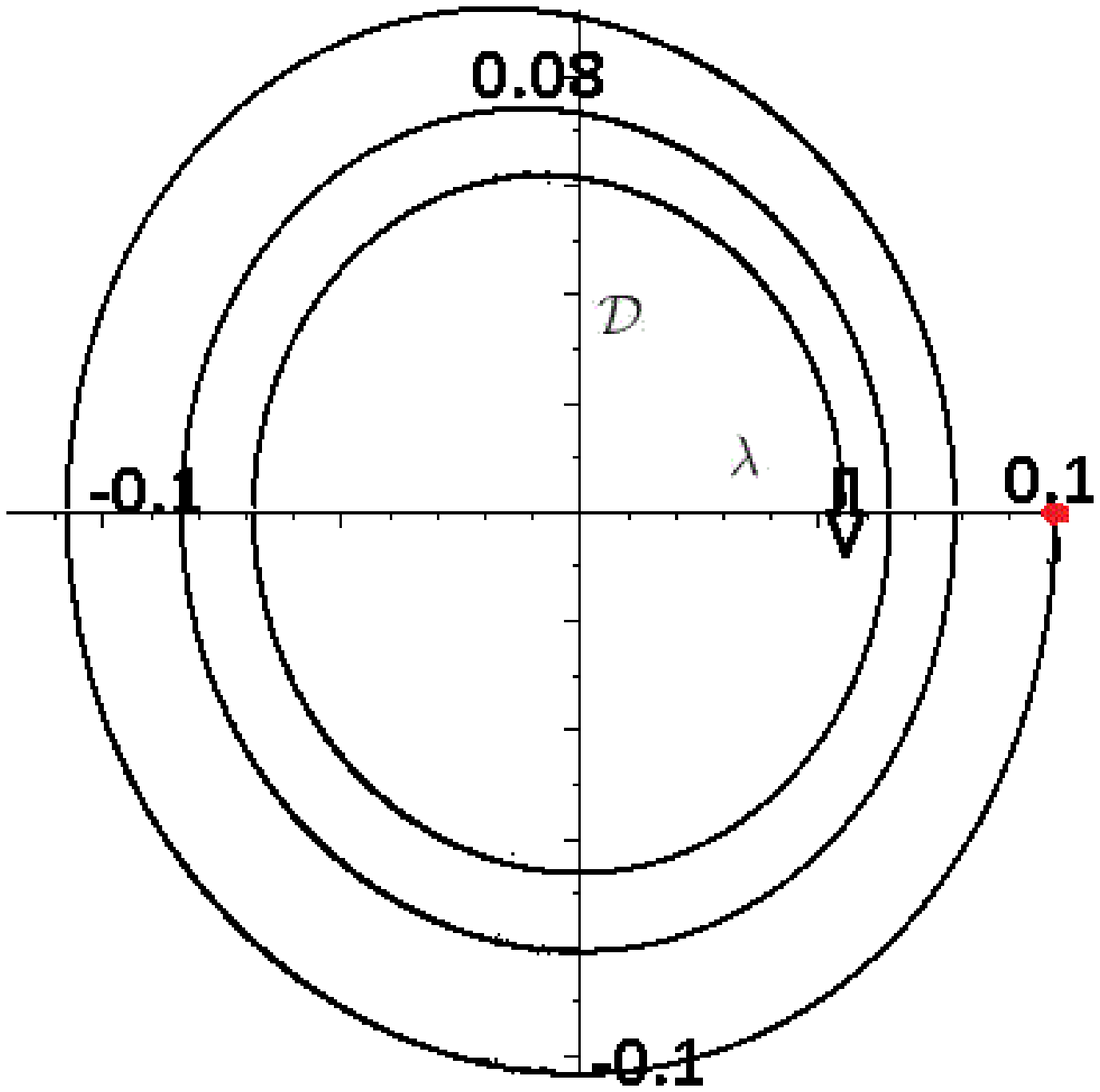}
\end{minipage}
\caption{Bounds for the projection of the phase trajectory on the
plane $(\lambda, \mathcal D)$. The starting point is $(0.1,0)$.
Left: the above bound $L$, the guaranteed number of oscillations is
$3$. Right: the respective below bound $l$.}\label{L}
\end{figure}

Fig.\ref{L}, left, presents the curve $L$, which bounds the
projection of the phase trajectory on the plane $(\lambda, \mathcal
D)$ from above. Fig.\ref{L}, right, presents the curve $l$, which
bounds the projection of the phase trajectory on the plane
$(\lambda, \mathcal D)$ from below. The method guarantees 3
oscillations before the blow-up.  The guaranteed time of smoothness
$$T_*>\inf\limits_{r_0\in \mathbb R} T_l(r_0)=T_l(0)=18.8685...$$ The
estimate of this guaranteed time from above is $19.1298...$

\begin{remark}
A series of computations with the data \eqref{gauss} was performed
in \cite{GFCA}, \cite{CH18}, so we have an opportunity to compare
our results with  the results of  sophisticated numerics. In
\cite{GFCA} computations are made for $a_*=0.365$, $\rho_*=0.6$,
i.e. $K\approx 0.222$, the breaking time is about $\theta_*=35$
(dimensionless units). Thus, the solution keeps smoothness within
the first oscillation. However, Proposition \ref{prop}
 does not guarantee the preservation of smoothness during the first oscillation.
 This means that the sufficient condition from Proposition \ref{prop} is still too rough.
\end{remark}
\begin{remark}
The method for refining estimates outlined in this section can be
easily adapted to the 3D radially symmetric case.
\end{remark}

\section{Discussion}\label{SecD}

In this paper, we constructed a method, that allows one to obtain a sufficient condition that guarantees the boundedness of the density component on a given time interval for  the solution to the Cauchy problem of the system
of PDE describing general 3D electrostatic oscillations.
 For particular cases
(plain and irrotational oscollations) we find an explicit sufficient condition for preventing a blow-up in the first period of oscillations.
 Next, we compare the sufficient condition with the known criteria for the formation of singularities, if such a possibility exists. Further, we consider the case of
2D axisymmetric oscillations with special initial data, where there
is a possibility of comparison with the results of a numerical study
of the blow-up process.
For this case, we obtain a sufficient condition for the preservation of smoothness, which is less rough than in the general case, as well as a sufficient condition for a blow-up.

As shown in §\ref{Sec3}, in the 1D case, as well as for the affine and radially symmetric cases in many dimensions, they exist globally in $t$ smooth solutions of the Cauchy problem. It is very interesting to check whether there are globally smooth solutions in the multidimensional case without these restrictive assumptions.  The numerics suggest that any other solution
necessarily blows up. Recently, a finite-time blow-up for general radially symmetric initial data, even arbitrarily small, has been proved analytically\cite{R_PhysicaD22}. The only exception is the initial data in the form of so-called simple waves, where the solution either tends to the affine one or blows up. The present paper can be regarded as an addition to this result, since it makes it possible to estimate the lifetime of a smooth solution. Now it is not known where affine solutions without radial symmetry necessarily blows up.

Another important open question is the description of the complete class of solutions corresponding to electrostatic oscillations.
The solutions considered in this article belong to this class, and the hypothesis is that in the multidimensional case, electrostatic oscillations must necessarily be radially symmetric or affine. However, this fact has not yet been proven.

\section*{Acknowledgment}
Supported by  the Moscow Center for
Fundamental and Applied Mathematics under the agreement
№075-15-2019-1621.

\end{document}